\documentclass[9pt,twoside,reqno]{amsart}
\usepackage{amssymb}
\usepackage{amsmath}
\usepackage{amsfonts}

\newtheorem{theorem}{Theorem}
\theoremstyle{plain}

\newtheorem{definition}{Definition}
\newtheorem{example}{Example}
\newtheorem{lemma}{Lemma}

\newtheorem{remark}{Remark}

\numberwithin{equation}{section}
\input{tcilatex}

\begin{document}
\title[]{Convergence rate of implicit iteration process and a data
dependence result}
\author{Isa Yildirim}
\address{Department of Mathematics, Ataturk University, Erzurum 25240, Turkey%
}
\email{isayildirim@atauni.edu.tr}
\author{Mujahid Abbas}
\address{Department of Mathematics, Lahore University of Management
Sciences, 54792 Lahore, Pakistan.}
\email{mujahid@lums.edu.pk}

\begin{abstract}
The aim of this paper is to introduce an implicit S-iteration process and
study its convergence in the framework of W-hyperbolic spaces. We show that
the implicit S-iteration process has higher rate of convergence than
implicit Mann type iteration and implicit Ishikawa-type iteration processes.
We present a numerical example to support the analytic result proved herein.
Finally, we prove a data dependence result for a contractive type mapping
using implicit S-iteration process.
\end{abstract}

\maketitle

\section{\textbf{Introduction and preliminaries}}

Throughout this paper, the letter $%
\mathbb{N}
$ will denote the set of natural numbers.

The theory of fixed points deals with the conditions which guarantee that a
mapping $T$ of \ a set $X$ into itself admits one or more fixed points, that
is, points $x$ of $X$ which solve an operator equation $x=Tx,$ called a
fixed point equation$.$ Fixed point theory serves as an essential tool for
solving problems arising in various branches of mathematical analysis. Over
the past two decades the development of fixed point theory in metric spaces
has attracted considerable attention due to numerous applications in areas
such as variational and linear inequalities, optimization, and approximation
theory.

The set $\left\{ p\in X:p=Tp\right\} $ of all fixed points of $T$ is denoted
by $F(T).$

One of the basic and the most widely applied fixed point theorem in all of
analysis is "Banach ( or Banach- Cassioppoli ) Contraction principle" \ due
to Banach \cite{Banach}. This principle lies at the heart of metric fixed
point theory. It states that if $(X,d)$ is a complete metric space and $%
T:X\rightarrow X$ satisfies 
\begin{equation*}
d(Tx,Ty)\leq kd(x,y),
\end{equation*}%
for all $x,y\in X,$ with $k\in (0,1),$ then $T$ has a unique fixed point. \
Due to its applications in mathematics and other related disciplines, Banach
contraction principle has been generalized in many directions.

Zamfirescu \cite{a37} obtained an important generalization of Banach fixed
point theorem using Zamfirescu mapping.

A self mapping $T$ on a metric space $X$ is called Zamfirescu mapping if
there exist real numbers $a,b,c$ satisfying $0<a<1,0<b,c<1/2$ such that at
least one of the following is true:%
\begin{equation}
\left\{ 
\begin{array}{cc}
\left( z_{1}\right) & d\left( Tx,Ty\right) \leq ad\left( x,y\right) , \\ 
\left( z_{2}\right) & d\left( Tx,Ty\right) \leq b\left( d\left( x,Tx\right)
+d\left( y,Ty\right) \right) , \\ 
\left( z_{3}\right) & d\left( Tx,Ty\right) \leq c\left( d\left( x,Ty\right)
+d\left( y,Tx\right) \right) .%
\end{array}%
\right.  \label{1}
\end{equation}

for any $x,y\in X.$

The contractive condition (\ref{1}) can be reformulated as follows:%
\begin{equation}
\left\{ 
\begin{array}{cc}
\left( b_{1}\right) & d\left( Tx,Ty\right) \leq \delta d\left( x,y\right)
+2\delta d\left( x,Tx\right) \text{ if one uses }\left( z_{2}\right) ,\text{
and } \\ 
\left( b_{2}\right) & d\left( Tx,Ty\right) \leq \delta d\left( x,y\right)
+2\delta d\left( x,Ty\right) \text{ if one uses }\left( z_{3}\right) ,%
\end{array}%
\right.  \label{2}
\end{equation}%
for all $x,y\in X,$ where $\delta =\max \left\{ a,\frac{b}{1-b},\frac{c}{1-c}%
\right\} $ (\cite{a13})$.$ Clearly, $\delta \in \left[ 0,1\right) .$

A mapping satisfying condition $\left( b_{1}\right) $ or $\left(
b_{2}\right) $ is called a quasi-contractive mapping. This class of mappings
is general than the class of Zamfirescu mappings.

Osilike and Udomene \cite{a20} extended the above class of mappings and
introduced a mapping $T$ satisfying the following contractive condition:%
\begin{equation}
\text{\ }d\left( Tx,Ty\right) \leq \delta d\left( x,y\right) +Ld\left(
x,Tx\right) ,  \label{3}
\end{equation}%
for all $x,y\in X,$ where $L\geq 0$ and $\delta \in \left[ 0,1\right) .$

For more results and discussion in this direction, we refer to \cite{abbas}, 
\cite{khan-common} and references mentioned therein.

Imoru and Olantiwo \cite{a38} gave the following definition:

\begin{definition}
A self mapping $T$ on $X$ is called a contractive-like mapping if there
exists a constant $\delta \in \left[ 0,1\right) $ and a strictly increasing
and continuous function $\varphi :\left[ 0,\infty \right) \rightarrow \left[
0,\infty \right) $ with $\varphi \left( 0\right) =0$ such that for any $%
x,y\in X,$ we have%
\begin{equation}
d\left( Tx,Ty\right) \leq \delta d\left( x,y\right) +\varphi \left( d\left(
x,Tx\right) \right) .  \label{4}
\end{equation}
\end{definition}

\begin{definition}
\cite{soltuz} Let $T,S$ be two self mappings on $X$. We say that $S$ is an
approximate operator of $T$ if for all $\varepsilon >0$, we have $d\left(
Tx,Sx\right) \leq \varepsilon $ holds for any $x\in X$ .
\end{definition}

An ambient space equipped with certain convexity structure play a
significant role in solving a fixed point equation. Since Banach space is a
vector space, the concept of a line segment joining any two points of a
nonempty subset of a Banach space give rise to the convexity structure. \
However, metric spaces do not naturally have this convex structure. The
notion of convex metric spaces was introduced by Takahashi \cite{Ta} who
studied the fixed points of nonexpansive mappings in the setting of such
spaces. All normed spaces and their convex subsets are convex metric spaces.
But there are many examples of convex metric spaces which are not embedded
in any normed space (\cite{Ta}). Over time, different convex structures have
been introduced on metric spaces.

Kohlenbach \cite{koh} introduced $W-$hyperbolic spaces as follows:

\begin{definition}
A $W-$hyperbolic space $(X,d,W)$ is a metric space $(X,d)$ together with a
convexity mopping $W:$ $X^{2}\times \lbrack 0,1]\rightarrow X$ satisfying
the following properties:%
\begin{equation*}
\begin{array}{l}
\text{(i) }d(u,W(x,y,\alpha ))\leq (1-\alpha )d(u,x)+\alpha d(u,y) \\ 
\text{(ii) }d(W(x,y,\alpha ),W(x,y,\beta ))=\left\vert \alpha -\beta
\right\vert d(x,y) \\ 
\text{(iii) }W(x,y,\alpha )=W(y,x,1-\alpha ) \\ 
\text{(iv) }d(W(x,z,\alpha ),W(y,w,\alpha ))\leq (1-\alpha )d(x,y)+\alpha
d(z,w)%
\end{array}%
\end{equation*}%
for all $x,y,z,w\in X$ and $\alpha ,\beta \in \lbrack 0,1].$
\end{definition}

If the triplet $(X,d,W)$ satisfies the condition (i) only, then it coincides
with the convex metric space introduced by Takahashi \cite{Ta}. Every
hyperbolic space is a convex metric space but converse does not hold in
general ( \cite{a1}). A subset $E$ of a $W-$hyperbolic space $X$ is convex
if $W(x,y,\alpha )\in E$ for all $x,y\in E$ and $\alpha \in \lbrack 0,1].$%
\newline
Note that every $W-$hyperbolic space is a geodesic space. CAT(0) spaces,
normed linear space, The Hilbert ball and Busseman spaces are important
examples of $W-$ hyperbolic spaces. A $W-$ hyperbolic space represents a
unified approach for both linear and nonlinear structures simultaneously.
There are hyperbolic spaces which are not imbedded in any Banach space .

On the other hand, different iterative algorithms have been used to
approximate the solution of a fixed point equation (\cite{Chang, AR, Kim,
Liu, Ta}). Implicit iterative schemes are of great importance from numerical
stand point as they provide accurate approximation ( see, \cite{a1}, \cite%
{Khan}, \cite{cir, cir1, yil}).

The motivation of this paper is to define an implicit S-iteration process
with higher rate of convergence when compared with Mann type(\ref{ek3}) and
Ishikawa type (\ref{ek2}) implicit iterative processes.

Let $E$ be a nonempty convex subset of a $W-$hyperbolic space $X$ and $%
T:E\rightarrow E$. $\ $Choose $x_{0}\in E$ and define the sequence $%
\{x_{n}\} $ as follows:%
\begin{eqnarray}
x_{n} &=&W\left( Tx_{n-1},Ty_{n},\alpha _{n}\right)  \label{de0} \\
y_{n} &=&W\left( x_{n},Tx_{n},\beta _{n}\right) ,\text{ \ }n\in 
\mathbb{N}
,  \notag
\end{eqnarray}%
where $\left\{ a_{n}\right\} $ and $\left\{ \beta _{n}\right\} $ are certain
real sequences in $\left[ 0,1\right] $.

We translate implicit Ishikawa and implicit Mann iteration processes
introduced by \'{C}iri\'{c} et al. (\cite{cir1, cir3}) in the setup of $W-$%
hyperbolic space as follows: 
\begin{eqnarray}
x_{n} &=&W\left( x_{n-1},Ty_{n},\alpha _{n}\right)  \label{ek2} \\
y_{n} &=&W\left( x_{n},Tx_{n},\beta _{n}\right) \text{ \ }n\in 
\mathbb{N}
,  \notag
\end{eqnarray}%
and%
\begin{equation}
x_{n}=W\left( x_{n-1},Ty_{n},\alpha _{n}\right) ,\text{ }n\in 
\mathbb{N}
.  \label{ek3}
\end{equation}

\begin{remark}
Note that, the process (\ref{de0}) is independent of (\ref{ek2}) and (\ref%
{ek3}) in the sense that neither of them reduce to the other.
\end{remark}

The following definition is due to Berinde \cite{ber}.

\begin{definition}
Let $\left\{ x_{n}\right\} $ and $\left\{ u_{n}\right\} $ be two fixed point
iteration processes in a metric space $X$ such that $\lim\limits_{n%
\rightarrow \infty }$ $x_{n}=\lim\limits_{n\rightarrow \infty } $ $u_{n}=p,$
where $p$ is a fixed point of a self mapping $T$ on $X.$ Suppose that%
\begin{equation*}
d(x_{n},p)\leq a_{n}\text{ and }d(u_{n},p)\leq b_{n}\text{, }n\in 
\mathbb{N}
.
\end{equation*}
where $\left\{ a_{n}\right\} $ and $\left\{ b_{n}\right\} $ are two null
sequences of positive numbers. If $\left\{ a_{n}\right\} $ converges faster
than $\left\{ b_{n}\right\} $, then we say $\left\{ x_{n}\right\} $
converges faster than $\left\{ u_{n}\right\} $ to $p$.
\end{definition}

We also need the following lemma in order to prove our main results.

\begin{lemma}
\label{l1}\cite{soltuz} Let ${\left\{ a_{n}\right\} }$ be a nonnegative
sequence. If there exists an $n_{0}\in 
\mathbb{N}
$ such that for all $n\geq n_{0}$, we have 
\begin{equation*}
a_{n+1}\leq (1-\mu _{n})a_{n}+\mu _{n}\eta _{n},
\end{equation*}
where $\mu _{n}\in \left( 0,1\right) $, $\sum_{n=0}^{\infty }\mu _{n}=\infty 
$ and $\eta _{n}\geq 0$ for all $n\in 
\mathbb{N}
$. Then the following holds:%
\begin{equation*}
0\leq \lim_{n\rightarrow \infty }\sup a_{n}\leq \lim_{n\rightarrow \infty
}\sup \eta _{n}.
\end{equation*}
\end{lemma}

\section{\textbf{Main Results}}

We start with the following result.

\begin{theorem}
\label{main1}Let $E$ be a nonempty closed convex subset of $W-$hyperbolic
space $X$ and $T:E\rightarrow E$ a contractive type mapping with $F\left(
T\right) \neq \emptyset $. Then, for the sequence $\{x_{n}\}$ defined in ($%
\ref{de0}$) with $\sum \left( 1-\alpha _{n}\right) =\infty $, we have $%
\lim\limits_{n\rightarrow \infty }x_{n}=p,$ where $p\in F\left( T\right) $.

\begin{proof}
Suppose that $p\in F\left( T\right) $. Using ($\ref{de0}$) and ($\ref{4}$),
we have\ \ \ \ \ \ \ \ \ \ \ \ \ \ \ \ \ \ \ \ \ \ \ \ \ \ \ \ \ \ \ \ \ \ \
\ \ \ \ \ \ \ \ \ \ \ \ \ \ \ \ \ \ \ \ \ \ \ \ \ \ \ \ \ \ \ \ \ \ \ \ \ \
\ \ \ \ \ \ \ \ \ \ \ \ \ \ \ \ \ \ \ \ \ \ \ \ \ \ \ \ \ \ \ \ \ \ \ \ \ \
\ \ \ \ \ \ \ \ \ \ \ \ \ \ \ \ \ \ \ \ \ \ \ \ \ \ \ \ \ \ \ \ \ \ \ \ \ \
\ \ \ \ \ \ \ \ \ \ \ \ \ \ \ \ \ \ \ \ \ \ \ \ \ \ \ \ \ \ \ \ \ \ \ \ \ \
\ \ \ \ \ \ \ \ \ \ \ \ \ \ \ \ \ \ \ \ \ \ \ \ \ \ \ \ \ \ \ \ \ 
\begin{eqnarray}
d(x_{n},p) &=&d(W\left( Tx_{n-1},Ty_{n},\alpha _{n}\right) ,p)  \label{tr1}
\\
&\leq &\alpha _{n}d(Tx_{n-1},p)+(1-\alpha _{n})d(Ty_{n},p)  \notag \\
&\leq &\alpha _{n}\left[ \delta d(x_{n-1},p)+\varphi \left( d\left(
p,Tp\right) \right) \right]  \notag \\
&&+(1-\alpha _{n})\left[ \delta d(y_{n},p)+\varphi \left( d\left(
p,Tp\right) \right) \right]  \notag \\
&=&\alpha _{n}\delta d(x_{n-1},p)+(1-\alpha _{n})\delta d(y_{n},p)  \notag
\end{eqnarray}

and%
\begin{eqnarray}
d(y_{n},p) &=&d(W\left( x_{n},Tx_{n},\beta _{n}\right) ,p)  \label{tr2} \\
&\leq &\beta _{n}d(x_{n},p)+(1-\beta _{n})d(Tx_{n},p)  \notag \\
&\leq &\beta _{n}d(x_{n},p)+(1-\beta _{n})\left[ \delta d(x_{n},p)+\varphi
\left( d\left( p,Tp\right) \right) \right]  \notag \\
&=&\beta _{n}d(x_{n},p)+(1-\beta _{n})\delta d(x_{n},p)  \notag \\
&=&\left[ \beta _{n}+(1-\beta _{n})\delta \right] d(x_{n},p).  \notag
\end{eqnarray}%
Therefore,%
\begin{equation*}
d(x_{n},p)\leq \alpha _{n}\delta d(x_{n-1},p)+(1-\alpha _{n})\delta \left[
\beta _{n}+(1-\beta _{n})\delta \right] d(x_{n},p).
\end{equation*}%
That is, 
\begin{equation}
\left[ 1-(1-\alpha _{n})\delta \left[ \beta _{n}+(1-\beta _{n})\delta \right]
\right] d(x_{n},p)\leq \alpha _{n}\delta d(x_{n-1},p)\text{,}  \label{tre4}
\end{equation}%
which further implies that 
\begin{equation}
d(x_{n},p)\leq \frac{\alpha _{n}\delta }{1-(1-\alpha _{n})\delta \left[
\beta _{n}+(1-\beta _{n})\delta \right] }d(x_{n-1},p).  \label{denk3}
\end{equation}%
We set%
\begin{equation*}
\Delta _{n}=\frac{\alpha _{n}\delta }{1-(1-\alpha _{n})\delta \left[ \beta
_{n}+(1-\beta _{n})\delta \right] }.
\end{equation*}%
Then%
\begin{eqnarray*}
1-\Delta _{n} &=&1-\frac{\alpha _{n}\delta }{1-(1-\alpha _{n})\delta \left[
\beta _{n}+(1-\beta _{n})\delta \right] } \\
&=&\frac{1-(1-\alpha _{n})\delta \left[ \beta _{n}+(1-\beta _{n})\delta %
\right] -\alpha _{n}\delta }{1-(1-\alpha _{n})\delta \left[ \beta
_{n}+(1-\beta _{n})\delta \right] } \\
&\geq &1-(1-\alpha _{n})\delta \left[ \beta _{n}+(1-\beta _{n})\delta \right]
-\alpha _{n}\delta
\end{eqnarray*}
implies that%
\begin{eqnarray}
\Delta _{n} &\leq &(1-\alpha _{n})\delta \left[ \beta _{n}+(1-\beta
_{n})\delta \right] +\alpha _{n}\delta  \label{tr5} \\
&=&(1-\alpha _{n})\delta \left[ 1-\left( 1-\delta \right) (1-\beta _{n})%
\right] +\alpha _{n}\delta  \notag \\
&\leq &(1-\alpha _{n})\delta +\alpha _{n}  \notag \\
&=&1-(1-\alpha _{n})\left( 1-\delta \right) .  \notag
\end{eqnarray}%
By ($\ref{denk3}$) and ($\ref{tr5}$)$,$ we have%
\begin{eqnarray}
d(x_{n},p) &\leq &\left[ 1-(1-\alpha _{n})\left( 1-\delta \right) \right]
d(x_{n-1},p)  \label{tr6} \\
&\leq &\dprod\limits_{i=1}^{n}\left[ 1-(1-\alpha _{i})\left( 1-\delta
\right) \right] d(x_{0},p).  \notag
\end{eqnarray}%
Since $a>0,$ $1+a\leq e^{a},$ $(\ref{tr6})$ gives that 
\begin{eqnarray}
d(x_{n},p) &\leq &\exp \left\{ \tsum\limits_{i=1}^{n}(1-\alpha _{i})\left(
1-\delta \right) \right\} d(x_{0},p)  \label{tr7} \\
&\leq &\exp \left\{ \tsum\limits_{n=1}^{\infty }(1-\alpha _{n})\left(
1-\delta \right) \right\} d(x_{0},p).  \notag
\end{eqnarray}%
Using the fact that $0\leq \delta <1$ and $\dsum (1-\alpha _{n})=\infty $,
we conclude that $\lim_{n\rightarrow \infty }d(x_{n},p)=0$.
\end{proof}
\end{theorem}

The following result deals with the rate of convergence of implicit
S-iteration process.

\begin{theorem}
Let $E$ be a nonempty closed convex subset of $W-$hyperbolic space $X$ and $%
T:E\rightarrow E$ a contractive type mapping with $F\left( T\right) \neq
\emptyset $. Then, the sequence $\{x_{n}\}$ defined in ($\ref{de0}$) with $%
\sum \left( 1-\alpha _{n}\right) =\infty $ converges to the fixed point of $%
T $ faster than implicit Ishikawa type (\ref{ek2}) and implicit Mann type (%
\ref{ek3}) iterations.

\begin{proof}
Let $p$ be a fixed point of $T$. Using an implicit Ishikawa type (\ref{ek2})
iteration process we have,%
\begin{eqnarray}
d(x_{n},p) &=&d(W\left( x_{n-1},Ty_{n},\alpha _{n}\right) ,p)  \label{tr8} \\
&\leq &\alpha _{n}d(x_{n-1},p)+(1-\alpha _{n})d(Ty_{n},p)  \notag \\
&\leq &\alpha _{n}d(x_{n-1},p)+(1-\alpha _{n})\left[ \delta
d(y_{n},p)+\varphi \left( d\left( p,Tp\right) \right) \right]  \notag \\
&=&\alpha _{n}d(x_{n-1},p)+(1-\alpha _{n})\delta d(y_{n},p)  \notag
\end{eqnarray}%
and%
\begin{eqnarray}
d(y_{n},p) &=&d(W\left( x_{n},Tx_{n},\beta _{n}\right) ,p)  \label{tr10} \\
&\leq &\beta _{n}d(x_{n},p)+(1-\beta _{n})d(Tx_{n},p)  \notag \\
&\leq &\beta _{n}d(x_{n},p)+(1-\beta _{n})\left[ \delta d(x_{n},p)+\varphi
\left( d\left( p,Tp\right) \right) \right]  \notag \\
&=&\beta _{n}d(x_{n},p)+(1-\beta _{n})\delta d(x_{n},p)  \notag \\
&=&\left[ \beta _{n}+\delta (1-\beta _{n})\right] d(x_{n},p).  \notag
\end{eqnarray}%
Therefore, 
\begin{equation}
d(x_{n},p)\leq \alpha _{n}d(x_{n-1},p)+(1-\alpha _{n})\delta \left[ \beta
_{n}+\delta (1-\beta _{n})\right] d(x_{n},p).  \label{tr11}
\end{equation}%
That is,%
\begin{eqnarray}
d(x_{n},p) &\leq &\frac{\alpha _{n}}{1-(1-\alpha _{n})\delta \left[ \beta
_{n}+\delta (1-\beta _{n})\right] }d(x_{n-1},p)  \label{tr12} \\
&\leq &...\leq c_{n}  \notag
\end{eqnarray}%
where%
\begin{equation}
c_{n}=\left( \frac{\alpha _{n}}{1-(1-\alpha _{n})\delta \left[ \beta
_{n}+\delta (1-\beta _{n})\right] }\right) ^{n}d(x_{0},p).  \label{tr13}
\end{equation}%
Using implicit Mann iteration (\ref{ek2}), we obtain that%
\begin{eqnarray}
d(x_{n},p) &=&d(W\left( x_{n-1},Tx_{n},\alpha _{n}\right) ,p)  \label{tre13}
\\
&\leq &\alpha _{n}d(x_{n-1},p)+(1-\alpha _{n})d(Tx_{n},p)  \notag \\
&\leq &\alpha _{n}d(x_{n-1},p)+(1-\alpha _{n})\left[ \delta
d(x_{n},p)+\varphi \left( d\left( p,Tp\right) \right) \right]  \notag \\
&=&\alpha _{n}d(x_{n-1},p)+(1-\alpha _{n})\delta d(x_{n},p).  \notag
\end{eqnarray}
Thus, we have%
\begin{eqnarray}
d(x_{n},p) &\leq &\frac{\alpha _{n}}{1-(1-\alpha _{n})\delta }d(x_{n-1},p)
\label{tr131} \\
&\leq &b_{n}  \notag
\end{eqnarray}%
where%
\begin{equation}
b_{n}=\left( \frac{\alpha _{n}}{1-(1-\alpha _{n})\delta }\right)
^{n}d(x_{0},p).  \label{tr132}
\end{equation}%
Using implicit S-iteration process and following the arguments of the proof
of Theorem 1, we have 
\begin{eqnarray*}
d(x_{n},p) &\leq &\frac{\alpha _{n}\delta }{1-(1-\alpha _{n})\delta \left[
\beta _{n}+(1-\beta _{n})\delta \right] \left[ \gamma _{n}+(1-\gamma
_{n})\delta \right] }d(x_{n-1},p) \\
&\leq &...\leq a_{n}
\end{eqnarray*}%
where%
\begin{equation}
a_{n}=\left( \frac{\alpha _{n}\delta }{1-(1-\alpha _{n})\delta \left[ \beta
_{n}+(1-\beta _{n})\delta \right] \left[ \gamma _{n}+(1-\gamma _{n})\delta %
\right] }\right) ^{n}d(x_{0},p).  \label{tr14}
\end{equation}%
Note that $\lim_{n\rightarrow \infty }\frac{a_{n}}{c_{n}}=0$ and $%
\lim_{n\rightarrow \infty }\frac{a_{n}}{b_{n}}=0$. Indeed, 
\begin{eqnarray*}
\left[ \gamma _{n}+(1-\gamma _{n})\delta \right] &<&1\Rightarrow \left[
\beta _{n}+(1-\beta _{n})\delta \right] <1 \\
&\Rightarrow &(1-\alpha _{n})\delta \left[ \beta _{n}+(1-\beta _{n})\delta %
\right] <1 \\
&\Rightarrow &1-(1-\alpha _{n})\delta \left[ \beta _{n}+(1-\beta _{n})\delta %
\right] >0,
\end{eqnarray*}%
give

\begin{eqnarray*}
\alpha _{n}\delta &<&\alpha _{n}\Rightarrow \frac{\alpha _{n}\delta }{%
1-(1-\alpha _{n})\delta \left[ \beta _{n}+\delta (1-\beta _{n})\right] } \\
&<&\frac{\alpha _{n}}{1-(1-\alpha _{n})\delta } \\
&\Rightarrow &\left( \frac{\alpha _{n}\delta }{1-(1-\alpha _{n})\delta \left[
\beta _{n}+\delta (1-\beta _{n})\right] }\right) ^{n} \\
&<&\left( \frac{\alpha _{n}}{1-(1-\alpha _{n})\delta }\right) ^{n}
\end{eqnarray*}%
and%
\begin{eqnarray*}
\alpha _{n}\delta &<&\alpha _{n}\Rightarrow \frac{\alpha _{n}\delta }{%
1-(1-\alpha _{n})\delta \left[ \beta _{n}+\delta (1-\beta _{n})\right] } \\
&<&\frac{\alpha _{n}}{1-(1-\alpha _{n})\delta \left[ \beta _{n}+\delta
(1-\beta _{n})\right] } \\
&\Rightarrow &\left( \frac{\alpha _{n}\delta }{1-(1-\alpha _{n})\delta \left[
\beta _{n}+\delta (1-\beta _{n})\right] }\right) ^{n} \\
&<&\left( \frac{\alpha _{n}}{1-(1-\alpha _{n})\delta \left[ \beta
_{n}+\delta (1-\beta _{n})\right] }\right) ^{n}.
\end{eqnarray*}%
Hence, we have%
\begin{eqnarray*}
\lim_{n\rightarrow \infty }\frac{a_{n}}{c_{n}} &=&\lim_{n\rightarrow \infty }%
\frac{\left( \frac{\alpha _{n}\delta }{1-(1-\alpha _{n})\delta \left[ \beta
_{n}+\delta (1-\beta _{n})\right] }\right) ^{n}d(x_{0},p)}{\left( \frac{%
\alpha _{n}}{1-(1-\alpha _{n})\delta }\right) ^{n}d(x_{0},p)} \\
&=&0,
\end{eqnarray*}%
and%
\begin{eqnarray*}
\lim_{n\rightarrow \infty }\frac{a_{n}}{b_{n}} &=&\lim_{n\rightarrow \infty }%
\frac{\left( \frac{\alpha _{n}\delta }{1-(1-\alpha _{n})\delta \left[ \beta
_{n}+\delta (1-\beta _{n})\right] }\right) ^{n}d(x_{0},p)}{\left( \frac{%
\alpha _{n}}{1-(1-\alpha _{n})\delta \left[ \beta _{n}+\delta (1-\beta _{n})%
\right] }\right) ^{n}d(x_{0},p)} \\
&=&0.
\end{eqnarray*}
\end{proof}
\end{theorem}

We now support our above analytical proof by a numerical example using
MATLAB.

\begin{example}
\label{example6}Let $E=\left[ 0,1\right] $ and $T:E\rightarrow E$ a mapping
defined by $Tx=\frac{x}{2}.$ Note that $T$ is a contractive type mapping$.$%
Choose $\alpha _{n}=1-\frac{1}{n}$ and $\beta _{n}=1-\frac{1}{n},n\geq 2$
and for $n=1,$ $\alpha _{n}=\beta _{n}=0.$ The comparison of the
convergences of the implicit S-iteration ($\ref{de0}$), implicit Ishikawa
type (\ref{ek2}) and implicit Mann type iterations (\ref{ek3}) to the fixed
point $p=0$ are given in table with the initial value $x_{1}=1$.

\ The following table presents a\textbf{\ }comparison of rate of convergence
of the implicit S-iteration process with implicit Ishikawa type and implicit
Mann type iteration processes for the mapping \ given in Example \ref%
{example6}.%
\begin{equation*}
\begin{tabular}{|c|c|c|c|}
\hline
\textbf{n}\  & \textbf{IMI} & \textbf{III} & \textbf{ISI} \\ \hline
2 & 0.666666666666667 & 0.615384615384615 & 0.307692307692308 \\ \hline
5 & 0.406349206349206 & 0.352704628530670 & 0.022044039283167 \\ \hline
7 & 0.340992340992341 & 0.292145335107371 & 0.004564770861053 \\ \hline
10 & 0.283773192751521 & 0.240691952056443 & 0.000470101468860 \\ \hline
13 & 0.248169351176485 & 0.209336831746067 & 0.000051107624938 \\ \hline
16 & 0.223294138742407 & 0.187699995568689 & 0.000005728149279 \\ \hline
20 & 0.199408653447441 & 0.167113839554526 & 0.000000318744353 \\ \hline
25 & 0.178133771931084 & 0.148920204678483 & 0.000000008876336 \\ \hline
30 & 0.162477710197415 & 0.135609685643003 & 0.000000000252593 \\ \hline
35 & 0.150335628473559 & 0.125328510781087 & 0.000000000007295 \\ \hline
40 & 0.140563343828096 & 0.117078595772533 & 0.000000000000213 \\ \hline
43 & 0.135541774913220 & 0.112847389889567 & 0.000000000000026 \\ \hline
46 & 0.131022580805197 & 0.109043978938918 & 0.000000000000003 \\ \hline
50 & 0.125645129018549 & 0.104523598655989 & 0.000000000000000 \\ \hline
\end{tabular}%
\end{equation*}

\begin{remark}
From the example above, we see that our iteration ISI (implicit S-iteration)
is faster than the III (implicit Ishikawa iteration) and IMI (implicit Mann
iteration) under the same control conditions.
\end{remark}
\end{example}

Finally, we present a data dependence result for contractive type mapping $T$
using implicit S-iteration process.

\begin{theorem}
\label{main1 copy(1)}Let $E$ be a nonempty closed convex subset of $W-$%
hyperbolic space $X$, $T:E\rightarrow E$ a contractive type mapping, $S$ an
approximate operator of $T$ and $\{x_{n}\}$ a sequence defined in ($\ref{de0}
$)$.$ For contractive type mapping $T$, define an iterative sequence $%
\{u_{n}\}$ as follows:%
\begin{eqnarray}
u_{n} &=&W\left( Su_{n-1},Tv_{n},\alpha _{n}\right)  \label{ekson} \\
v_{n} &=&W\left( u_{n},Su_{n},\beta _{n}\right) ,\text{ \ }n\in 
\mathbb{N}
,  \notag
\end{eqnarray}%
where $\left\{ a_{n}\right\} $ and $\left\{ \beta _{n}\right\} $ are real
sequences in $\left[ 0,1\right] $ satisfying $a_{n}<1$, $n\in 
\mathbb{N}
$, and $\sum \left( 1-\alpha _{n}\right) =\infty $. If $Tp=p$ and $Sq=q$
such that $u_{n}\rightarrow q$ as $n\rightarrow \infty $, then we have, $%
d\left( p,q\right) \leq \dfrac{2\varepsilon }{\left( 1-\delta \right) ^{2}}$%
, $\varepsilon >0$ is some appropriate number.

\begin{proof}
By ($\ref{4}$), ($\ref{de0}$) and ($\ref{ekson}$), we have\ \ \ \ \ \ \ \ \
\ \ \ \ \ \ \ \ \ \ \ \ \ \ \ \ \ \ \ \ \ \ \ \ \ \ \ \ \ \ \ \ \ \ \ \ \ \
\ \ \ \ \ \ \ \ \ \ \ \ \ \ \ \ \ \ \ \ \ \ \ \ \ \ \ \ \ \ \ \ \ \ \ \ \ \
\ \ \ \ \ \ \ \ \ \ \ \ \ \ \ \ \ \ \ \ \ \ \ \ \ \ \ \ \ \ \ \ \ \ \ \ \ \
\ \ \ \ \ \ \ \ \ \ \ \ \ \ \ \ \ \ \ \ \ \ \ \ \ \ \ \ \ \ \ \ \ \ \ \ \ \
\ \ \ \ \ \ \ \ \ \ \ \ \ \ \ \ \ \ \ \ \ \ \ \ \ \ \ \ \ \ \ \ \ \ \ \ \ \
\ \ \ \ \ \ \ \ \ \ \ \ \ \ \ \ \ \ \ \ \ 
\begin{eqnarray}
d(y_{n},v_{n}) &=&d\left( W\left( x_{n},Tx_{n},\beta _{n}\right) ,W\left(
u_{n},Su_{n},\beta _{n}\right) \right)  \label{ekson1} \\
&\leq &\beta _{n}d\left( x_{n},u_{n}\right) +(1-\beta _{n})d\left(
Tx_{n},Su_{n}\right)  \notag \\
&\leq &\beta _{n}d\left( x_{n},u_{n}\right) +(1-\beta _{n})\left[ d\left(
Tx_{n},Su_{n}\right) +d\left( Tu_{n},Su_{n}\right) \right]  \notag \\
&\leq &\beta _{n}d\left( x_{n},u_{n}\right) +(1-\beta _{n})\left[ \delta
d\left( x_{n},u_{n}\right) +\varphi \left( d\left( x_{n},Tx_{n}\right)
\right) +\varepsilon \right]  \notag \\
&=&\left[ \beta _{n}+(1-\beta _{n})\delta \right] d\left( x_{n},u_{n}\right)
+(1-\beta _{n})\left[ \varphi \left( d\left( x_{n},Tx_{n}\right) \right)
+\varepsilon \right]  \notag
\end{eqnarray}%
and%
\begin{eqnarray}
d(x_{n},u_{n}) &=&d\left( W\left( Tx_{n-1},Ty_{n},\alpha _{n}\right)
,W\left( Su_{n-1},Sv_{n},\alpha _{n}\right) \right)  \label{ekson2} \\
&\leq &\alpha _{n}d\left( Tx_{n-1},Su_{n-1}\right) +(1-\alpha _{n})d\left(
Ty_{n},Sv_{n}\right)  \notag \\
&\leq &\alpha _{n}\left[ d\left( Tx_{n-1},Tu_{n-1}\right) +d\left(
Tu_{n-1},Su_{n-1}\right) \right]  \notag \\
&&+(1-\alpha _{n})\left[ d\left( Ty_{n},Tv_{n}\right) +d\left(
Tv_{n},Sv_{n}\right) \right]  \notag \\
&\leq &\alpha _{n}\left[ \delta d\left( x_{n-1},u_{n-1}\right) +\varphi
\left( d\left( x_{n-1},Tx_{n-1}\right) \right) +\varepsilon \right]  \notag
\\
&&+(1-\alpha _{n})\left[ \delta d\left( y_{n},v_{n}\right) +\varphi \left(
d\left( y_{n},Ty_{n}\right) \right) +\varepsilon \right]  \notag \\
&=&\alpha _{n}\delta d\left( x_{n-1},u_{n-1}\right) +\alpha _{n}\varphi
\left( d\left( x_{n-1},Tx_{n-1}\right) \right) +\varepsilon  \notag \\
&&+(1-\alpha _{n})\delta d\left( y_{n},v_{n}\right) +(1-\alpha _{n})\varphi
\left( d\left( y_{n},Ty_{n}\right) \right) .  \notag
\end{eqnarray}%
Substituting ($\ref{ekson1}$) in ($\ref{ekson2}$), we obtain that%
\begin{eqnarray}
d(x_{n},u_{n}) &\leq &\alpha _{n}\delta d\left( x_{n-1},u_{n-1}\right)
+\alpha _{n}\varphi \left( d\left( x_{n-1},Tx_{n-1}\right) \right)
\label{ekson3} \\
&&+(1-\alpha _{n})\varphi \left( d\left( y_{n},Ty_{n}\right) \right)
+\varepsilon +(1-\alpha _{n})\delta  \notag \\
&&\left[ \left[ \beta _{n}+(1-\beta _{n})\delta \right] d\left(
x_{n},u_{n}\right) +(1-\beta _{n})\left[ \varphi \left( d\left(
x_{n},Tx_{n}\right) \right) +\varepsilon \right] \right]  \notag \\
&=&\alpha _{n}\delta d\left( x_{n-1},u_{n-1}\right) +\alpha _{n}\varphi
\left( d\left( x_{n-1},Tx_{n-1}\right) \right)  \notag \\
&&+(1-\alpha _{n})\varphi \left( d\left( y_{n},Ty_{n}\right) \right)
+\varepsilon +(1-\alpha _{n})\delta  \notag \\
&&\left[ \beta _{n}+(1-\beta _{n})\delta \right] d\left( x_{n},u_{n}\right)
+(1-\alpha _{n})\delta  \notag \\
&&(1-\beta _{n})\varphi \left( d\left( x_{n},Tx_{n}\right) \right)
+(1-\alpha _{n})\delta (1-\beta _{n})\varepsilon .  \notag
\end{eqnarray}%
It follows from ( $\ref{ekson3}$) that%
\begin{eqnarray}
&&\left[ 1-(1-\alpha _{n})\delta \left[ \beta _{n}+(1-\beta _{n})\delta %
\right] \right] d(x_{n},u_{n})  \label{ekson4} \\
&\leq &\alpha _{n}\delta d\left( x_{n-1},u_{n-1}\right) +\alpha _{n}\varphi
\left( d\left( x_{n-1},Tx_{n-1}\right) \right)  \notag \\
&&+(1-\alpha _{n})\varphi \left( d\left( y_{n},Ty_{n}\right) \right)
+\varepsilon +(1-\alpha _{n})\delta  \notag \\
&&(1-\beta _{n})\varphi \left( d\left( x_{n},Tx_{n}\right) \right)
+(1-\alpha _{n})\delta (1-\beta _{n})\varepsilon ,  \notag
\end{eqnarray}%
which further implies that 
\begin{eqnarray}
&&d(x_{n},u_{n})\leq \frac{\alpha _{n}\delta }{1-(1-\alpha _{n})\delta \left[
\beta _{n}+(1-\beta _{n})\delta \right] }d\left( x_{n-1},u_{n-1}\right)
\label{ekson5} \\
&&+\frac{\left\{ 
\begin{array}{c}
\alpha _{n}\varphi \left( d\left( x_{n-1},Tx_{n-1}\right) \right) +(1-\alpha
_{n})\varphi \left( d\left( y_{n},Ty_{n}\right) \right) \\ 
+(1-\alpha _{n})\delta (1-\beta _{n})\varphi \left( d\left(
x_{n},Tx_{n}\right) \right)%
\end{array}%
\right\} }{1-(1-\alpha _{n})\delta \left[ \beta _{n}+(1-\beta _{n})\delta %
\right] }  \notag \\
&&\frac{+\varepsilon +(1-\alpha _{n})\delta (1-\beta _{n})\varepsilon }{%
1-(1-\alpha _{n})\delta \left[ \beta _{n}+(1-\beta _{n})\delta \right] }. 
\notag
\end{eqnarray}%
We set%
\begin{equation*}
\Delta _{n}=\frac{\alpha _{n}\delta }{1-(1-\alpha _{n})\delta \left[ \beta
_{n}+(1-\beta _{n})\delta \right] }.
\end{equation*}%
Following arguments similar to those in the proof of Theorem $\ref{main1}$,
we have%
\begin{equation}
\Delta _{n}\leq 1-(1-\alpha _{n})\left( 1-\delta \right) .  \label{ekson6}
\end{equation}%
Therefore,%
\begin{eqnarray*}
d(x_{n},u_{n}) &\leq &\left[ 1-(1-\alpha _{n})\left( 1-\delta \right) \right]
d\left( x_{n-1},u_{n-1}\right) \\
&&+\frac{(1-\alpha _{n})\left( 1-\delta \right) \left\{ 
\begin{array}{c}
\frac{\alpha _{n}}{1-\alpha _{n}}\varphi \left( d\left(
x_{n-1},Tx_{n-1}\right) \right) +\varphi \left( d\left( y_{n},Ty_{n}\right)
\right) \\ 
+\delta (1-\beta _{n})\varphi \left( d\left( x_{n},Tx_{n}\right) \right)
+2\varepsilon%
\end{array}%
\right\} }{\left( 1-\delta \right) \left[ 1-(1-\alpha _{n})\delta \left[
\beta _{n}+(1-\beta _{n})\delta \right] \right] }.
\end{eqnarray*}%
Note that%
\begin{eqnarray*}
1-(1-\alpha _{n})\delta \left[ \beta _{n}+(1-\beta _{n})\delta \right]
&=&1-(1-\alpha _{n})\delta \left[ 1-(1-\beta _{n})\left( 1-\delta \right) %
\right] \\
&\geq &1-\delta .
\end{eqnarray*}%
Hence%
\begin{eqnarray}
d(x_{n},u_{n}) &\leq &\left[ 1-(1-\alpha _{n})\left( 1-\delta \right) \right]
d\left( x_{n-1},u_{n-1}\right)  \label{ekson7} \\
&&+\frac{(1-\alpha _{n})\left( 1-\delta \right) \left\{ 
\begin{array}{c}
\frac{\alpha _{n}}{1-\alpha _{n}}\varphi \left( d\left(
x_{n-1},Tx_{n-1}\right) \right) +\varphi \left( d\left( y_{n},Ty_{n}\right)
\right) \\ 
+\delta (1-\beta _{n})\varphi \left( d\left( x_{n},Tx_{n}\right) \right)
+2\varepsilon%
\end{array}%
\right\} }{\left( 1-\delta \right) ^{2}}.  \notag
\end{eqnarray}%
That is,%
\begin{equation*}
a_{n+1}\leq (1-\mu _{n})a_{n}+\mu _{n}\eta _{n}
\end{equation*}%
where $a_{n+1}=d(x_{n},u_{n})$, $\mu _{n}=(1-\alpha _{n})\left( 1-\delta
\right) $ and%
\begin{equation*}
\eta _{n}=\frac{%
\begin{array}{c}
\frac{\alpha _{n}}{1-\alpha _{n}}\varphi \left( d\left(
x_{n-1},Tx_{n-1}\right) \right) +\varphi \left( d\left( y_{n},Ty_{n}\right)
\right) \\ 
+\delta (1-\beta _{n})\varphi \left( d\left( x_{n},Tx_{n}\right) \right)
+2\varepsilon%
\end{array}%
}{\left( 1-\delta \right) ^{2}}.
\end{equation*}

From Theorem $\ref{main1}$, we have $\lim_{n\rightarrow \infty
}d(x_{n},p)=\lim_{n\rightarrow \infty }d(x_{n-1},p)=0$, and $%
\lim_{n\rightarrow \infty }d(u_{n},p)=0$. As $\varphi $ is continuous, $%
\lim_{n\rightarrow \infty }\varphi \left( d\left( x_{n-1},Tx_{n-1}\right)
\right) =0$, $\lim_{n\rightarrow \infty }\varphi \left( d\left(
y_{n},Ty_{n}\right) \right) =0$, and $\lim_{n\rightarrow \infty }\varphi
\left( d\left( x_{n},Tx_{n}\right) \right) =0$. Thus all the conditions of
Lemma $\ref{l1}$ are satisfied, therefore ($\ref{ekson7}$) becomes%
\begin{equation*}
d\left( p,q\right) \leq \frac{2\varepsilon }{\left( 1-\delta \right) ^{2}}.
\end{equation*}
\end{proof}
\end{theorem}


\begin{thebibliography}{99}
\bibitem{abbas} M. Abbas, P. Vetro, S. H. Khan, On fixed points of Berinde's
contractive mappings in cone metric spaces, Carpath. J. Math., 26(2) (2010)
121-133.

\bibitem{Banach} S. Banach, Sur les op\'{e}rations dans les ensembles
abstraits et leur applications aux \'{e}quations int\'{e}grales, Fund.
Math., 3 (1922), 133-181.

\bibitem{ber} V. Berinde, Picard iteration converges faster than Mann
iteration for a class of quasicontractive operators, Fixed Point Theory and
Applications, 2004(2004), 97-105.

\bibitem{a13} V. Berinde, On the convergence of the Ishikawa iteration in
the class of quasi contractive operators. Acta Math. Univ. Comen. 73,
119-126 (2004).

\bibitem{Chang} S.S. Chang, L. Yang, X.R. Wang, Stronger convergence
theorems for an infinite family of uniformly quasi-Lipschitzian mappings in
convex metric spaces, Appl. Math. Comp. 217 (2010) 277--282.

\bibitem{a1} R. Chugh, P. Malik and V. Kumar, On analytical and numerical
study of implicit fixed point iterations, Cogent Mathematics (2015), 2:
1021623.

\bibitem{cir} Lj. B. \'{C}iri\'{c}, Rafiq, A., Caki\'{c}, N., \& Ume, J. S.
(2009). Implicit Mann fixed point iterations for pseudo-contractive
mappings. Applied Mathematics Letters, 22, 581--584.

\bibitem{cir1} Lj. B. \'{C}iri\'{c}, Rafiq, A., Radenovi\'{c}, S., Rajovi%
\'{c}, M., \& Ume, J. S. (2008). On Mann implicit iterations for strongly
accretive and strongly pseudo-contractive mappings. Applied Mathematics and
Computation, 198, 128--137.

\bibitem{cir3} Lj. B.Ciric, Ume, J. S. M., \& Khan, S. (2003). On the
convergence of the Ishikawa iterates to a common fixed point of two
mappings. Archivum Mathematicum (Brno) Tomus, 39, 123--127.

\bibitem{a38} C.O. Imoru, M.O. Olantiwo, On the stability of Picard and Mann
iteration processes. Carpath. J. Math. 19, 155-160 (2003).

\bibitem{AR} A.R. Khan, M.A. Ahmed, Convergence of a general iterative
scheme for a finite family of asymptotically quasi-nonexpansive mappings in
convex metric spaces and applications, Comput. Math. Appl. 59 (2010)
2990-2995.

\bibitem{Khan} S.H. Khan, I. Yildirim, M. Ozdemir. Convergence of an
implicit algorithm for two families of nonexpansive mappings, Comput. Math.
Appl. 59 (2010) 3084-3091.{}

\bibitem{khan-common} S. H. Khan, Common fixed points of quasi-contractive
type operators by a generalized iterative process, IAENG Int. J. Appl.
Math., 41(3) (2011), 260-264.

\bibitem{Kim} J. K. Kim, K. S. Kim, S. M. Kim, Convergence theorems of
implicit iteration process for for finite family of asymptotically
quasi-nonexpansive mappings in convex metric space, Nonlinear Analysis and
Convex Analysis, 1484 (2006) 40-51.

\bibitem{koh} U. Kohlenbach, (2004). Some logical metatherems with
applications in functional analysis. Transactions of the American
Mathematical Society, 357, 89--128.

\bibitem{Liu} Q.Y. Liu, Z.B. Liu, N.J. Huang, Approximating the common fixed
points of two sequences of uniformly quasi-Lipschitzian mappings in convex
metric spaces, Appl. Math. Comp. 216 (2010) 883--889.

\bibitem{a20} M.O. Osilike, A. Udomene, Short proofs of stability results
for fixed point iteration procedures for a class of contractive-type
mappings. Indian J. Pure Appl. Math. 30, 1229-1234 (1999).

\bibitem{soltuz} S.M. \c{S}oltuz, Grosan, T., Data dependence for Ishikawa
iteration when dealing with contractive like operators. Fixed Point Theory
Appl. 2008, Article ID 242916 (2008). doi:10.1155/2008/242916

\bibitem{Ta} W. Takahashi, A convexity in metric space and nonexpansive
mappings, Kodai Math. Sem. Rep. 22 (1970) 142-149.

\bibitem{yil} I. Yildirim, S. H. Khan, \textquotedblleft Convergence
theorems for common fixed points of asymptotically quasi-nonexpansive
mappings in convex metric spaces\textquotedblright , Applied Mathematics and
Computation, Volume 218, Issue 9 (2012), 4860-4866.

\bibitem{a37} T. Zamfirescu, Fix point theorems in metric spaces. Arch.
Math. 23, 292-298 (1972).
\end{thebibliography}
\end{document}